\documentclass[a4paper,11pt]{article}

\usepackage{MyLaTeX,MyArticle}
\input xy
\xyoption{all}
\renewcommand{\central}[2]{\mathrm{Z}_{#1}(#2)}
\title{Normal Subsystems of Fusion Systems}
\author{David A.~Craven, University of Oxford}
\date{June 2010}

\renewcommand{\F}{\mathcal{F}}
\newcommand{\E}{\mathcal{E}}

\begin{document}
\maketitle

\begin{abstract} In this article we prove that for any saturated fusion system, that the (unique) smallest weakly normal subsystem of it on a given strongly closed subgroup is actually normal. This has a variety of corollaries, such as the statement that the notion of a simple fusion system is independent of whether one uses weakly normal or normal subsystems. We also develop a theory of weakly normal maps, consider intersections and products of weakly normal subsystems, and the hypercentre of a fusion system.
\end{abstract}

\section{Introduction}

The theory of fusion systems is becoming an important topic in algebra, with interactions with group theory, representation theory and topology. This article is concerned with the structure of normal subsystems of fusion systems. There are two notions of a `normal' subsystem in the literature, one stronger than the other. (We will recall their definitions in this article.) We will follow \cite{craven2010} and call the subsystem considered by Aschbacher in \cite{aschbacher2008} a \emph{normal subsystem}, and the subsystem considered by, among others, Linckelmann in \cite{linckelmann2007} a \emph{weakly normal} subsystem.

Our first result highlights the exact relationship between normal and weakly normal subsystems.

\begin{thma}\label{thm:weaknormalnormal} Let $\F$ be a saturated fusion system on a finite $p$-group $P$. If $\E$ is a weakly normal subsystem of $\F$, then $\Orth^{p'}(\E)$ is a normal subsystem of $\F$.
\end{thma}

Recall that $\Orth^{p'}(\E)$ is the smallest (weakly) normal subsystem of $\E$ on the same $p$-group as $\E$ (see \cite[Section 6.5]{puig2006} or \cite[Theorem 5.4]{bcglo2007} for example, or Section \ref{sec:intersections} of this article). What Theorem \ref{thm:weaknormalnormal} says is that any weakly normal subsystem can be thought of as a normal subsystem $\E$, together with some $p'$-automorphisms of $\E$ that lie in $\F$. In other words, we have the following corollary.

\begin{cora} Let $\F$ be a saturated fusion system on a finite $p$-group $P$. If $\E$ is a weakly normal subsystem of $\F$, on a subgroup $T$, then there exists a normal subsystem $\E'$ of $\F$, also on $T$, and a $p'$-subgroup $H$ of $\Aut_\F(T)$ such that $\E$ is generated by $\E'$ and $H$.
\end{cora}

Theorem \ref{thm:weaknormalnormal} has several corollaries, such as the following.

\begin{cora} Let $\F$ be a saturated fusion system on a finite $p$-group $P$. Then $\F$ has no proper, non-trivial normal subsystems if and only if $\F$ has no proper, non-trivial weakly normal subsystems. In other words, the notion of a simple fusion system is the same with either normal subsystems or weakly normal subsystems.
\end{cora}

We give some more corollaries to Theorem \ref{thm:weaknormalnormal} in Section \ref{sec:cortothma}.

Let $\F$ be a saturated fusion system on a finite $p$-group $P$, and let $T$ be a strongly $\F$-closed subgroup of $P$. In \cite{aschbacher2008}, Aschbacher developed the concept of a normal map, which is a special type of function $A(-)$ on the set of all subgroups $U$ of $T$, with $A(U)$ a subgroup of $\Aut_\F(U)$. If $\E$ is a normal subsystem on $T$, then $A(U)=\Aut_\E(U)$ for $U\leq T$ is a normal map. Conversely, if $A(-)$ is a normal map, then the subsystem generated by $A(U)$ for $U\leq T$ is a normal subsystem $\E'$, and $A(U)=\Aut_{\E'}(U)$ for all $U\leq T$. In the proof of this theorem in \cite{aschbacher2008}, the fact that the map was associated to a normal subsystem was pivotal to the proof that the subsystem generated by the map was saturated. Here we decouple the theorem from this requirement, proving a similar result for weakly normal subsystems.

\begin{thma}\label{thm:weaknormalmap} Let $\F$ be a saturated fusion system on a finite $p$-group $P$, and let $T$ be a strongly $\F$-closed subgroup of $P$. Let $A(-)$ be a function from the set of subgroups $U$ of $T$, such that $A(U)$ is a subgroup of $\Aut_\F(U)$, satisfying the following conditions.
\begin{enumerate}
\item If $\phi$ is an $\F$-isomorphism whose domain is $U$, then $A(U\phi)=A(U)^\phi$.
\item If $U$ is fully $\F$-normalized, then $\Aut_T(U)\leq A(U)$.
\item $\Aut_T(T)$ is a Sylow $p$-subgroup of $A(T)$.
\item If $U$ is fully $\F$-normalized, then every element of $A(U)$ extends to an element of $A(U\Cent_T(U))$.
\item If $U$ is fully $\F$-normalized and $\Cent_T(U)\leq U$ then for any subgroup $U\leq V\leq \Norm_T(U)$, denoting by $A(U\leq V)$ the set of automorphisms of $V$ that restrict to automorphisms of $U$, the restriction map 
\[ A(U\leq V)\to \Norm_{A(U)}(\Aut_V(U)) \]
is surjective.
\end{enumerate}
If $\E$ is the subsystem generated by $A(U)$ for all $U\leq T$, then $\E$ is a weakly normal subsystem of $\F$ and $\Aut_\E(U)=A(U)$ for all $U\leq T$.
\end{thma}

When amalgamated with the axiom needed for weakly normal subsystems to be normal, we get a (slightly) weaker formulation of the axioms for $A(-)$ to be a normal map.

Let $\F$ be a saturated fusion system on a finite $p$-group $P$. The centre $\centre\F$ of $\F$ seems to play an important role in the theory of fusion systems: for example, in \cite[Corollary 6.14]{bcglo2007} it is proved that $\F$ is a fusion system of a finite group if and only if $\F/\centre\F$ is, and in \cite[Lemma 8.10]{aschbacher2008} it is proved that there is a one-to-one correspondence between normal subsystems of $\F$ on subgroups containing $\centre\F$ and normal subsystems of $\F/\centre\F$. As with finite groups, write $\central1\F=\centre\F$ and $\central i\F$ for the preimage in $P$ of $\centre{\F/\central{i-1}\F}$. The series $(\central i\F)$ eventually stabilizes; write $\central\infty\F$ for this limit, called the \emph{hypercentre} of $\F$.

Another subsystem that seems important is $\Norm_P(Q)\Cent_\F(Q)$; if $Q$ is a fully $\F$-normalized subgroup of $P$, then the subsystem $\Norm_P(Q)\Cent_\F(Q)$ is the (saturated) subsystem of $\F$ on $\Norm_P(Q)$ consisting of all morphisms $\phi:R\to S$ in $\F$ such that $\phi$ extends to $\bar\phi:QR\to QS$ with $\bar\phi|_Q=c_g$ for some $g\in\Norm_P(Q)$, where $c_g$ denotes conjugation by $g$. This notion is connected to the hypercentre by the following theorem.

\begin{thma}\label{thm:hypercentral} Let $\F$ be a saturated fusion system on a finite $p$-group $P$. If $X$ is a normal subgroup of $P$, then $\F=P\Cent_\F(X)$ if and only if $X\leq \central\infty\F$.
\end{thma}

In Section \ref{sec:hypercentral} we prove a few results about the hypercentre and central extensions; for example, we prove the following result. (A saturated fusion system $\F$ is \emph{perfect} if there is no surjective morphism $\F\to\F_A(A)$ for any non-trivial abelian $p$-group $A$.)

\begin{propa} Let $\F$ be a saturated fusion system on a finite $p$-group $P$. If $\F$ is perfect then $\central 2\F=\centre\F$. 
\end{propa}

This proposition allows us to define the universal central extension of a perfect fusion system, exactly analogously to perfect finite groups.

We also extend (slightly) Glauberman's $Z^*$-theorem, rewriting it in the language of fusion systems at the same time.

\begin{propa} Let $G$ be a finite group with $\Orth_{p'}(G)=1$, and let $P$ be a Sylow $p$-subgroup of $G$. Writing $\F=\F_P(G)$, we have that $\central i\F=\central iG$.
\end{propa}

We end with three sections applying our results. The first of these deals with intersections and products of subsystems, the second of these deals with corollaries to Theorem \ref{thm:weaknormalnormal} and the final section gives a comparison of weakly normal and normal subsystems, giving two situations in which normal subsystems behave better than weakly normal subsystems.

\medskip

The notation used in this article is quickly becoming standard in this field, and we refer to \cite{craven2010} for many of the notational issues. We will define weakly normal subsystems and normal subsystems in the next section, as this is still relatively new terminology, but leave the standard definitions to \cite{craven2010} (see also \cite{linckelmann2007} and \cite{blo2003}, although note that in the former `fusion systems' there are referred to as `saturated fusion systems' here). Note that in this article maps and homomorphisms are composed from left to right.

\section{Definitions and Notation}

The definitions of fusion system and saturated fusion system, and fully normalized, fully centralized, centric, radical and essential subgroups, are as in \cite{craven2010}, and these are now standard in the literature. (Notice that our definitions of fully normalized and fully centralized subgroups differ from that of \cite{puig2006}.) The definitions of a fully automized subgroup and a receptive subgroup are not common yet, and so we give them now.

\begin{defn} Let $\F$ be a fusion system on a finite $p$-group $P$. A subgroup $Q$ of $P$ is \emph{fully $\F$-automized} if $\Aut_P(Q)$ is a Sylow $p$-subgroup of $\Aut_\F(Q)$. A subgroup $R$ is \emph{receptive} if, whenever $\phi:S\to R$ is an isomorphism in $\F$, $\phi$ extends to a map $\bar\phi:N_\phi\to P$, where $N_\phi$ is the (full) preimage of the subgroup $\Aut_P(S)\cap\Aut_P(R)^{\phi^{-1}}$ under the natural map $\Norm_P(S)\to\Aut(S)$.
\end{defn}

If $Q$ is a subgroup of $P$ then by $\Aut_P(Q)$ we mean the set of automorphisms of $Q$ induced by conjugation by the elements of $P$ (actually $\Norm_P(Q)$). We introduce the notation $\Aut(Q\leq R)$ for the set of all automorphisms of $R$ that restrict to an automorphism of $Q$, where $Q\leq R$. We also use the obvious extensions of notation $\Aut_\F(Q\leq R)$ and $\Aut_P(Q\leq R)$. A subgroup $Q$ is said to have the \emph{surjectivity property} if, for any subgroup $R$ with $Q\Cent_P(Q)\leq R\leq \Norm_P(Q)$, the map
\[ \Aut_\F(Q\leq R)\to \Norm_{\Aut_\F(Q)}(\Aut_R(Q))\]
obtained by restriction is surjective.

The notion of an $\F$-invariant subsystem was introduced by Puig (see \cite{puig2006}, where it is called `normal'), and a saturated, $\F$-invariant subsystem is called \emph{weakly normal} here (and `normal' in \cite{linckelmann2007}). Normal subsystems are defined in \cite{aschbacher2008}, with the addition of one more axiom, which involves how the subsystem is embedded with respect to the centralizer of the underlying subgroup.

\begin{defn} Let $\F$ be a saturated fusion system on a finite $p$-group $P$, and let $T$ be a strongly $\F$-closed subgroup. A subsystem $\E\leq\F$, on $T$, is \emph{$\F$-invariant} if, whenever $Q\leq R$ are subgroups of $T$, $\phi:Q\to R$ is a morphism in $\E$ and $\psi:R\to T$ is a morphism in $\F$, the composition $\psi^{-1}\phi\psi:Q\psi\to R\psi$ is a morphism in $\E$. If $\E$ is saturated and $\F$-invariant then $\E$ is \emph{weakly normal} in $\F$.

If, in addition, whenever $\phi$ is an $\E$-automorphism of $T$, there is an $\F$-automorphism $\bar\phi\in\Aut_\F(T\Cent_P(T))$ such that $[\bar\phi,\Cent_P(T)]\leq \centre T$, then $\E$ is said to be \emph{normal} in $\F$.
\end{defn}

We also need to fix our notation for quotient subsystems and morphisms, taken from \cite{craven2010}.

\begin{defn} Let $\F$ be a saturated fusion system on a finite $p$-group $P$, and let $T$ be a strongly $\F$-closed subgroup of $P$. Let $\F/T$ denote the fusion system on $P/T$, whose morphisms $\Hom_{\F/T}(Q/T,R/T)$ (for $T\leq Q,R\leq P$) are all morphisms $\bar\phi$ induced from morphisms $\phi:Q\to R$ in $\F$ by taking quotients by $T$ (since $\phi$ acts as an automorphism of $T$).
\end{defn}

Let $\F$ be a fusion system on a finite $p$-group $P$, and let $Q$ be a subgroup of $P$. The subsystem $\Norm_P(Q)\Cent_\F(Q)$ is the subsystem on $\Norm_P(Q)$ consisting of all morphisms $\phi:R\to S$ such that $\phi$ extends to a map $\bar\phi:QR\to QS$ such that $\phi|_Q=c_g$ for some $g\in\Norm_P(Q)$, where $c_g$ denotes the conjugation map by $g$.

Recall that the \emph{centre} of $\F$, denoted $\centre\F$, is the set of all $x\in P$ such that any morphism $\phi:Q\to R$ in $\F$ has an extension $\bar\phi:Q\gen x\to R\gen x$ such that $x\bar\phi=x$. We may iterate this construction.

\begin{defn} Let $\F$ be a saturated fusion system on a finite $p$-group $P$. Define $\central1\F=\centre \F$ and $\central i\F$ inductively by $\central i\F$ being the preimage in $P$ of $\centre{\F/\central{i-1}\F}$. The largest term of this ascending sequence is the \emph{hypercentre} of $\F$, and is denoted by $\central\infty\F$.
\end{defn}

We introduce a couple more definitions concerning the generation and saturation of fusion systems.

\begin{defn} Let $\F$ be a fusion system on a finite $p$-group $P$. Let $\mc H$ be a set of subgroups of $P$, and let $\mc Q$ be an $\F$-conjugacy class of subgroups of $P$.
\begin{enumerate}
\item We say that $\F$ is \emph{$\mc H$-generated} if $\F$ is the smallest fusion system on $P$ containing all $\F$-morphisms between elements of $\mc H$.
\item We say that $\mc Q$ is \emph{saturated} if it contains a fully automized, receptive subgroup $Q$ of $P$.
\end{enumerate}
\end{defn}

Finally, we have a definition needed for studying hypercentres.

\begin{defn} Let $\F$ be a saturated fusion system on a finite $p$-group $P$. We say that $\F$ is \emph{perfect} if there is no morphism $\F\to\F_A(A)$ of fusion systems, for any abelian $p$-group $A$,.
\end{defn}

An equivalent way of defining being perfect is to say that $\Orth^p(\F)=\F$, where $\Orth^p(\F)$ is the hyperfocal subsystem (see \cite{bcglo2007}).

\section{Results from the Literature}

Having given all of the definitions, we turn to the theorems in the literature that we will need in the next sections. We start with a well-known result on $p'$-automorphisms of $p$-groups that we will use so often we will not comment.

\begin{lem}[{{\cite[Corollary 5.3.3]{gorenstein}}}] Let $P$ be a finite $p$-group, and let $Q$ be a normal subgroup of $P$. The subgroup of $\Aut(P)$ of automorphisms that act trivially on both $Q$ and $P/Q$ is a $p$-group.
\end{lem}

The next two propositions and lemma deal with proving saturation in a fusion system.

\begin{prop}[{{\cite[Lemma 2.4]{bcglo2005}}}]\label{prop:someofHgeneration1} Let $\F$ be a fusion system on a finite $p$-group $P$, and let $\mc H$ be a union of $\F$-conjugacy classes of subgroups of $P$. Suppose that $\F$ is $\mc H$-generated, and that every $\F$-conjugacy class in $\mc H$ is saturated. Let $\mc Q$ be an $\F$-conjugacy class of subgroups of $P$ such that all subgroups of $P$ properly containing members of $\mc Q$ lie in $\mc H$.
\begin{enumerate}
\item Any fully normalized subgroup $Q$ in $\mc Q$ is fully centralized, and for any other subgroup $\tilde Q\in \mc Q$, there is an $\F$-morphism $\phi:\Norm_P(\tilde Q)\to\Norm_P(Q)$ with $\tilde Q\phi=Q$.
\item For any fully normalized subgroup $Q$ in $\mc Q$, if $\{Q\}$ is saturated in $\Norm_\F(Q)$, then $\mc Q$ is saturated in $\F$.
\end{enumerate}
\end{prop}

\begin{prop}[{{\cite[Lemmas 2.4 and 2.5]{bcglo2005}}}]\label{prop:someofHgeneration2} Let $\F$ be a fusion system on a finite $p$-group $P$, and let $\mc H$ be a union of $\F$-conjugacy classes of subgroups of $P$. Suppose that $\F$ is $\mc H$-generated, and that every $\F$-conjugacy class in $\mc H$ is saturated. Let $\mc Q$ be an $\F$-conjugacy class of subgroups of $P$ such that all subgroups of $P$ of smaller index than those in $\mc Q$ lie in $\mc H$. Let $\mc H'$ denote the set of subgroups of $\Norm_P(Q)$ strictly containing $Q$.
\begin{enumerate}
\item Every $\Norm_\F(Q)$-conjugacy class of subgroups of $\mc H'$ is saturated, and every $\F$-automorphism of $Q$ is the composition of (restrictions of) $\Norm_\F(Q)$-morphisms between elements of $\mc H'$.
\item If either $Q$ is not $\F$-centric, or $\Out_P(Q)\cap \Orth_p(\Out_\F(Q))$ is non-trivial, then the $\Norm_\F(Q)$-conjugacy class $\{Q\}$ is saturated.
\end{enumerate}
\end{prop}

We need a few different equivalent conditions for a fusion system to be saturated.

\begin{thm}[Roberts--Shpectorov {{\cite[Theorems 5.2, 5.3]{robshp2009}}}]\label{thm:altrobshp} Let $\F$ be a fusion system on a finite $p$-group $P$. Then $\F$ is saturated if and only if $P$ is fully $\F$-automized and every $\F$-conjugacy class of subgroups of $P$ contains a fully $\F$-normalized, receptive member, and this is equivalent to every $\F$-conjugacy class of subgroups of $P$ containing a fully $\F$-automized, receptive member.
\end{thm}

\begin{thm}[Cf.\ Puig {{\cite[Corollary 2.14]{puig2006}}}]\label{thm:saturationequiv} Let $\F$ be a fusion system on a finite $p$-group $P$, and suppose that $P$ is fully automized. If, for each $\F$-conjugacy class, there is some subgroup $Q$ such that
\begin{enumerate}
\item if $R$ is $\F$-conjugate to $Q$ then there is a map $\Norm_P(R)\to \Norm_P(Q)$ that restricts to a map $R\to Q$, and
\item $Q$ has the surjectivity property,
\end{enumerate}
then $\F$ is saturated. Conversely, if $\F$ is saturated then any fully normalized subgroup $Q$ has the surjectivity property, and if $R$ is $\F$-conjugate to $Q$ then there is a map $\Norm_P(R)\to\Norm_P(Q)$ that restricts to a map $R\to Q$.
\end{thm}

The statement here is not quite what is proved in \cite{puig2006} -- there it is required that every fully normalized subgroup satisfies (i) and (ii) above -- but it is easily seen to be equivalent.

Combining this theorem with results of Broto, Castellana, Grodal, Levi and Oliver from \cite{bcglo2005}, we get the following result.

\begin{thm}\label{thm:normalmapprelim} Let $\F$ be a fusion system on a finite $p$-group $P$ with $P$ fully automized, and let $\mc H$ be the set of all $\F$-centric subgroups of $P$. If $\F$ is generated by $\{\Aut_\F(U):U \in \mc H\}$, and in every $\F$-conjugacy class of $\F$-centric subgroups there is a fully normalized subgroup with the surjectivity property, then $\F$ is saturated.
\end{thm}
\begin{pf} Let $\mc H'$ denote the union of all saturated conjugacy classes of $\F$, and suppose that $\mc Q$ is a conjugacy class of subgroups of smallest index subject to not being in $\mc H'$. Let $\E$ be the subsystem of $\F$ generated by all morphisms in $\Aut_\F(R)$, as $R$ runs over all subgroups of $P$ of larger order than the subgroups in $\mc Q$. Notice that $\Hom_\F(R,S)=\Hom_\E(R,S)$ if $R$ has larger order than the subgroups in $\mc Q$, and so such $\E$-conjugacy classes are identical to the $\F$-conjugacy classes, and also saturated. In addition, $\mc Q$ forms a single $\E$-conjugacy class; in particular, a subgroup in $\mc Q$ is $\F$-centric if and only if it is $\E$-centric.

If $\mc Q$ consists of subgroups that are not $\F$-centric, then by Propositions \ref{prop:someofHgeneration1}(ii) and \ref{prop:someofHgeneration2}(ii) $\mc Q$ is saturated in $\E$. Thus we may suppose that $\mc Q$ consists of $\F$-centric subgroups. Let $Q$ be a fully $\F$-normalized subgroup in $\mc Q$ with the surjectivity property, guaranteed by hypothesis. Proposition \ref{prop:someofHgeneration1}(i), applied to $\mc Q$ and $\E$, is exactly the second requirement for Theorem \ref{thm:saturationequiv}, and so $\mc Q$ is saturated, as required.
\end{pf}

We end with a couple of results on quotients, weakly normal subsystems and normalizers.

\begin{thm}[Puig {{\cite[Proposition 6.6]{puig2006}}}]\label{thm:normalifffrattini} Let $\F$ be a saturated fusion system on a finite $p$-group $P$, and let $\E$ be a saturated subsystem of $\F$, on the strongly $\F$-closed subgroup $T$ of $P$. The following are equivalent:
\begin{enumerate}
\item $\E$ is weakly normal in $\F$; and
\item every $\phi\in\Aut_\F(T)$ induces an automorphism of $\E$ (via $\psi\mapsto\psi^\phi$) and every $\F$-morphism $\psi:A\to B$ with $A,B\leq T$ may be written as $\psi=\alpha\beta$, where $\alpha\in\Aut_\F(T)$ and $\beta\in\Hom_\E(A\alpha,B)$.
\end{enumerate}
\end{thm}

\begin{prop}[Kessar--Linckelmann {{\cite[Proposition 3.4]{kessarlinckelmann2008}}}]\label{pcfqbijection} Let $\F$ be a saturated fusion system on a finite $p$-group $P$, and let $Q$ be a normal subgroup of $P$ such that $\F=P\Cent_\F(Q)$. If $\E$ is a saturated subsystem on a subgroup $R$ with $Q\leq R\leq P$, then $\E=\Norm_\F(R)$ if and only if $\E/Q=\Norm_{\F/Q}(R/Q)$.
\end{prop}

Combined with the next lemma, this yields a useful tool for finding subgroups $R$ such that $\F=\Norm_\F(R)$.

\begin{lem}\label{normalizersame} Let $\F$ be a saturated fusion system on a finite $p$-group $P$, and let $T$ be a strongly $\F$-closed subgroup of $P$. Let $Q$ be a subgroup of $P$ containing $T$. We have that $\Norm_\F(Q)/T=\Norm_{\F/T}(Q/T)$, and in particular $\Norm_P(Q)/T=\Norm_{P/T}(Q/T)$.
\end{lem}
\begin{pf} Clearly $Tg$ normalizes $Q/T$ if and only if $g$ normalizes $Q$, and so $\Norm_P(Q)/T=\Norm_{P/T}(Q/T)$.

If $\phi:R\to S$ is a morphism in $\Norm_\F(Q)$, then $\phi$ extends to a morphism $\psi:RQ\to SQ$ such that $Q\psi=Q$. If $\bar\phi$ denotes the image of $\phi$ in $\F/T$ then $\bar\phi:RT/T\to ST/T$ extends to $\bar\psi:RQ/T\to SQ/T$, with $(Q/T)\bar\psi=Q/T$. Hence $\Norm_\F(Q)/T\subs \Norm_{\F/T}(Q/T)$.

Conversely, let $\bar\phi:A/T\to B/T$ be a morphism in $\Norm_{\F/T}(Q/T)$, and let $\bar\psi:AQ/T\to BQ/T$ be an extension of $\bar\phi$ such that $(Q/T)\bar\psi=Q/T$. Let $\psi:QA\to QB$ be a morphism in $\F$ whose image in $\F/T$ is $\bar\psi$. Since $T$ is strongly $\F$-closed and $(Q/T)\bar\psi=Q/T$, we must have that $Q\psi=Q$. Similarly, as $(A/T)\bar\psi=B/T$, we must have that $A\psi=B$. Hence $\phi=\psi|_A$ is a morphism in $\Norm_\F(Q)$. It is also clear that the image of $\phi$ in $\F/T$ is $\bar\phi$, and so $\Norm_{\F/T}(Q/T)\subs \Norm_\F(Q)/T$, completing the proof.
\end{pf}

\section{Preliminary Results}

\begin{lem}\label{normalsubgroupssub} Let $\F$ be a saturated fusion system on a finite $p$-group $P$, and let $\E$ be a saturated subsystem of $\F$, on the subgroup $Q$ of $P$. We have that 
\[\Orth_p(\E)\geq \Orth_p(\F)\cap Q.\]
\end{lem}
\begin{pf} Let $R=\Orth_p(\F)$; by \cite[(3.7)]{aschbacherbirm} (see also \cite[Proposition 3.4]{craven2010}), a strongly $\F$-closed subgroup is contained in $\Orth_p(\F)$ if and only if it possesses a central series each of whose terms is strongly $\F$-closed. Hence $R$ possesses a central series
\[ 1=R_0\leq R_1\leq \cdots \leq R_d=R\]
such that each $R_i$ is strongly $\F$-closed. We claim that $Q_i=Q\cap R_i$ is strongly $\E$-closed; in this case,
\[ 1=Q_0\leq Q_1\leq \cdots \leq Q_d=Q\cap R\]
is a central series for $Q\cap R$, each of whose terms is strongly $\E$-closed, yielding that $Q\cap R\leq \Orth_p(\E)$ (by another application of \cite[(3.7)]{aschbacherbirm}), as required. It remains to show that $Q_i$ is strongly $\E$-closed; however, any morphism in $\E$ that originates inside $Q_i=Q\cap R_i$ must have image inside $Q$ since $\E$ is a subsystem on $Q$, and must also lie in $R_i$ since $R_i$ is strongly $\F$-closed, and so $Q_i$ is strongly $\E$-closed.
\end{pf}

\begin{lem} Let $\F$ be a saturated fusion system on a finite $p$-group $P$. If $X$ is a strongly $\F$-closed subgroup of $P$ then so is $X\Cent_P(X)$.
\end{lem}
\begin{pf} It is easy to see that $\Cent_P(X)$ is strongly $\Norm_\F(X)$-closed, and hence $X\Cent_P(X)/X$ is strongly $\Norm_\F(X)/X$-closed by \cite[Theorem 6.1]{craven2010}. However, since $\Norm_\F(X)/X=\F/X$, by definition of $\F/X$, we see that $X\Cent_P(X)/X$ is strongly $\F/X$-closed, so that $X\Cent_P(X)$ is strongly $\F$-closed, by \cite[Theorem 6.1]{craven2010} again.
\end{pf}

\begin{lem}\label{lem:intswithstrclosed} Let $\F$ be a saturated fusion system on a finite $p$-group $P$, and let $\E$ be a weakly normal subsystem of $\F$, on the subgroup $T$ of $P$. If $Q$ is an $\E$-centric subgroup of $T$, then $Q\Cent_P(Q)\cap T=Q$ and $Q\Cent_P(Q)\cap T\Cent_P(T)=Q\Cent_P(T)$.
\end{lem}
\begin{pf} Since $Q\leq T$, by the modular law $Q\Cent_P(Q)\cap T=Q(\Cent_P(Q)\cap T)=Q\centre Q=Q$. Taking the product of both sides of the first equation by $\Cent_P(T)$, and by the modular law again,
\[ Q\Cent_P(T)= (Q\Cent_P(Q)\cap T)\Cent_P(T)=Q\Cent_P(Q)\cap T\Cent_P(T),\]
as needed.
\end{pf}

We introduce a piece of notation at this point for a subquotient that will be important in our work to prove Theorem \ref{thm:weaknormalnormal}.

\begin{defn} Let $P$ be a finite $p$-group, and let $T$ be a normal subgroup of $P$. If $Q$ is a subgroup of $T$ such that $\Cent_T(Q)\leq Q$, write $Y_Q$ for the subquotient $\centre Q\Cent_P(T)/\centre Q$ of $P$.
\end{defn}

\begin{prop}\label{prop:extnstoQCP(T)} Let $\F$ be a saturated fusion system on a finite $p$-group $P$, and let $\E$ be a weakly normal subsystem of $\F$, on the subgroup $T$ of $P$. Let $Q$ and $R$ be $\E$-centric subgroups of $T$, and let $\phi:Q\to R$ be an isomorphism in $\F$. Let $\psi$ be an extension of $\phi$ to $Q\Cent_P(T)$ in $\F$.
\begin{enumerate}
\item $\im\psi=R\Cent_P(T)$.
\item $Q\psi=R$, $\centre Q\psi=\centre R$, and $(\centre Q\Cent_P(T))\psi=\centre R\Cent_P(T)$.
\item if $g\in \Norm_T(Q)$, then $c_g$ acts trivially on $\Cent_P(T)$, and so on $Y_Q$.
\end{enumerate}
\end{prop}
\begin{pf}Let $\bar Q=Q\Cent_P(T)$ and $\bar R=R\Cent_P(T)$. Since $\bar Q=Q\Cent_P(Q)\cap T\Cent_P(T)$, and $T\Cent_P(T)$ is strongly $\F$-closed, $\bar Q\psi\leq T\Cent_P(T)$. Also, since $\bar Q\leq Q\Cent_P(Q)$, $\bar Q\psi\leq R\Cent_P(R)$. Hence $\bar Q\psi \leq R\Cent_P(R)\cap T\Cent_P(T)=\bar R$. If $|\bar Q|=|\bar R|$, then we have proved (i); however, $Q\cap \Cent_P(T)=\centre Q$ and $R\cap \Cent_P(T)=\centre R$, so that indeed they have the same orders.

To prove (ii), clearly $Q\psi=R$, so $\centre Q\psi=\centre R$, and certainly $\Cent_{\bar Q}(Q)=\bar Q\cap \Cent_P(Q)=\centre Q\Cent_P(T)$, and $\Cent_{\bar Q}(Q)\psi=\Cent_{\bar R}(R)$.

We move on to (iii). The map $c_g$ is an automorphism of $Q$, and $c_g$ acts trivially on $\Cent_P(T)$. By (ii), $c_g$ acts as an automorphism of $Y_Q$, and this action is trivial since the action on $\Cent_P(T)$ is trivial.
\end{pf}

We will use the following easy lemma often to move between the various possible interpretations of the condition that $[\phi,\Cent_P(T)]$ is contained in various subgroups of $P$. 

\begin{lem}\label{lem:equivtorightaction} Let $\F$ be a saturated fusion system on a finite $p$-group $P$, and let $\E$ be a weakly normal subsystem of $\F$, on the subgroup $T$ of $P$. Let $Q$ and $R$ be $\E$-centric subgroups of $T$, and let $\phi:Q\Cent_P(T)\to R\Cent_P(T)$ be a morphism in $\F$. The following are equivalent:
\begin{enumerate}
\item $[\phi,\Cent_P(T)]\leq T$;
\item $[\phi,\Cent_P(T)]\leq \centre R$; and
\item for $x\in\Cent_P(T)$, we have $(\centre Q x)\phi=\centre R x$.
\end{enumerate}
\end{lem}
\begin{pf} Let $x$ be an element of $\Cent_P(T)$. Firstly, notice that $[\phi,x^{-1}]=(x\phi)x^{-1}$.

Assume (i). Since $x$ centralizes $T$ it centralizes $Q$, and so $x\phi\in\Cent_P(R)$. Also, $x$ centralizes $R$, so that $[\phi,x^{-1}]\in\Cent_P(R)\cap T=\centre R$, proving that (i) implies (ii). Clearly (ii) implies (i).

Next, consider the induced action of $\phi$ on $Y_Q$, which maps it to $Y_R$ by Proposition \ref{prop:extnstoQCP(T)}(ii). As $\centre Q\phi=\centre R$, the induced action is that $\centre Qx$ is sent to $(\centre Qx)\phi=\centre R(x\phi)$. Hence $\centre R(x\phi)=\centre Rx$ if and only if $(x\phi)x^{-1}\in\centre R$, proving that (ii) is equivalent to (iii), as needed.
\end{pf}

We end with the following trivial lemma, which has been used in various guises in many papers in this field.

\begin{lem}\label{lem:extnpauto} Let $U$ and $V$ be subgroups of the finite $p$-group $P$, and suppose that $U\normal V$. If $\phi$ is an automorphism of $V$ that acts trivially on $U$, and $\Cent_V(U)\leq U$, then $\phi$ is a $p$-automorphism.
\end{lem}
\begin{pf}Since $\phi$ acts trivially on $U$, the induced action on $\Aut(U)$ is also trivial, and hence $\phi$ acts trivially on $V/\centre U$, so also on $V/U$. The set of all automorphisms in $\Aut(V)$ that act trivially on both $U$ and $V/U$ is a $p$-group, and so in particular $\phi$ is a $p$-automorphism.
\end{pf}

\section{Weakly Normal and Normal Subsystems}

In this section, $\F$ is a saturated fusion system on a finite $p$-group $P$ and $T$ is a strongly $\F$-closed subgroup of $P$. There is a weakly normal subsystem $\E$ on $T$ with $\Orth^{p'}(\E)=\E$. Let $W$ denote the subset of $\Aut_\E(T)$ consisting of those automorphisms $\phi$ that extend to an automorphism $\bar\phi\in\Aut_\F(T\Cent_P(T))$ such that $[\bar\phi,\Cent_P(T)]\leq \centre T$.

Let $\alpha$ be an $\E$-automorphism of $T$. If $Q$ is a subgroup of $T$, call $Q$ a \emph{detecting subgroup} for $\alpha$ if there is a morphism $\beta:Q\Cent_P(T)\to R\Cent_P(T)$ (where $R=Q\alpha$) in $\F$ such that $\beta|_Q=\alpha|_Q$ and $[\beta,\Cent_P(T)]\leq \centre R$. Lemma \ref{lem:equivtorightaction} gives some equivalent conditions to $Q$ being a detecting subgroup for $\alpha$ in terms of the condition on $[\beta,\Cent_P(T)]$, and we will move between these conditions regularly.

Let $\E'$ denote the subsystem (not necessarily saturated) generated by $\Orth^{p'}(\Aut_\E(Q))$, as $Q$ runs over all subgroups of $T$. In \cite[Section 5]{bcglo2007}, Broto, Castellana, Grodal, Levi and Oliver prove that for every element $\alpha$ of $\Aut_\E(T)$, there exist $\E$-centric subgroups $Q$ and $R$ of $T$, with $\alpha|_Q:Q\to R$ in $\E'$. We will show that such a subgroup $Q$ is a detecting subgroup for $\alpha$.

\begin{prop}If $\alpha$ is an automorphism in $\Aut_\E(T)$ then there is a detecting subgroup for $\alpha$.
\end{prop}
\begin{pf} Notice that any morphism in $\E'$ is a composition of the restriction of $p$-automorphisms of subgroups of $T$, with these $p$-automorphisms lying in $\E$. We begin by examining such automorphisms.

Let $g$ be an element of $T$, and suppose that $g$ normalizes an $\E$-centric subgroup $S$ of $T$. Clearly $g$ centralizes $\Cent_P(T)$, and so induces an automorphism $c_g\in\Aut_\F(S\Cent_P(T))$.

If $\phi\in\Aut_\E(S)$ is any other automorphism, then $\phi$ extends to an $\F$-automorphism of $S\Cent_P(T)$. We claim that $(c_g)^\phi$ acts trivially on $Y_S=\centre S\Cent_P(T)/\centre S$. This is clear since $\phi$ acts on $Y_S$ by Proposition \ref{prop:extnstoQCP(T)}(ii), and $c_g$ acts trivially on $Y_S$ by (iii) of the proposition. Hence any $p$-automorphism of $S$ in $\E'$ extends to an $\F$-automorphism $\alpha$ of $S\Cent_P(T)$ acting trivially on $Y_S$. By Lemma \ref{lem:equivtorightaction}, this is equivalent to $[\alpha,\Cent_P(T)]\leq T$.

Let $\alpha\in\Aut_\E(T)$, and let $Q$ be an $\E$-centric subgroup such that $\alpha|_Q$ is a morphism in $\E'$. The isomorphism $\alpha|_Q:Q\to R$ is a composition of $p$-automorphisms $\theta_i$ of subgroups $S_i$ of $T$ in $\E'$. Since, for each $\theta_i$ there is an extension $\bar\theta_i\in\Aut_\F(S_i\Cent_P(T))$ with $[\bar\theta_i,\Cent_P(T)]\leq T$, taking the composition of these extensions and restricting to $Q\Cent_P(T)$ yields an isomorphism $\bar\alpha:Q\Cent_P(T)\to R\Cent_P(T)$ in $\F$ (whose image is such by Proposition \ref{prop:extnstoQCP(T)}(i)). We evaluate $[\bar\alpha,\Cent_P(T)]$: since each $[\bar\theta_i,\Cent_P(T)]$ lies inside $T$, and $[\bar\alpha,\Cent_P(T)]$ is contained in the product of conjugates of $[\bar\theta_i,\Cent_P(T)]$, we see that $[\bar\alpha,\Cent_P(T)]\leq T$, so that $[\bar\alpha,\Cent_P(T)]\leq \centre R$, proving that $Q$ is a detecting subgroup for $\alpha$.
\end{pf}

Now that we have proved that detecting subgroups exist, we choose $Q$ to be a detecting subgroup for an automorphism $\alpha$. The aim is to show that $\Norm_T(Q)$ is also a detecting subgroup for $\alpha$, for then $T$ is a detecting subgroup for $\alpha$, and so $\alpha$ lies in the subgroup $W$ defined at the beginning of the section. We will prove this result for $p'$-automorphisms of $T$ first, and then extend to all $\E$-automorphisms.

\begin{prop}\label{prop:caseRfullnormal} Let $\phi:Q\to R$ be an isomorphism in $\F$, where $Q$ and $R$ are $\E$-centric subgroups of $T$, and let $\psi$ be an extension of $\phi$ to $Q\Cent_P(T)$ in $\F$. Suppose that $N_\phi$ contains $\Norm_T(Q)$. If $R\Cent_P(T)$ is fully $\F$-normalized then $N_\psi$ contains $\Norm_T(Q)$, and so $\psi$ extends to $\theta:\Norm_T(Q)\Cent_P(T)\to \Norm_T(R)\Cent_P(T)$ in $\F$.
\end{prop}
\begin{pf} Let $g$ be an element of $\Norm_T(Q)$, and let $h\in\Norm_T(R)$ such that $c_g^\phi=c_h$. We claim that, as an element of $\Aut_\F(R\Cent_P(T))$, $c_g^\psi c_h^{-1}$ acts trivially on $R$ and $Y_R=\centre R\Cent_P(T)/\centre R$. To see this, obviously $c_g^\psi$ and $c_h$ act the same on $R$; also, $c_h$ acts trivially on $Y_R$ by Proposition \ref{prop:extnstoQCP(T)}(iii), so we need that $c_g^\psi$ acts trivially on $Y_R$. However, $c_g$ acts trivially on $Y_Q$ by Proposition \ref{prop:extnstoQCP(T)}(iii), and $\psi$ maps $Y_Q$ to $Y_R$ by Proposition \ref{prop:extnstoQCP(T)}(ii), so that this holds.

Let $K$ denote the set of all elements of $\Aut_\F(R\Cent_P(T))$ that act trivially on $R$ and $Y_R$. It is easy to see that $K$ is a $p$-group, because a $p'$-automorphism acting trivially on $\centre R$ and $\centre R\Cent_P(T)/\centre R$ acts trivially on $\centre R\Cent_P(T)$, so that (since it acts trivially on $R$) it acts trivially on $R\Cent_P(T)$. In fact, since every $\F$-automorphism of $R\Cent_P(T)$ induces an automorphism of $R$ and $Y_R$, $K$ is a normal $p$-subgroup of $\Aut_\F(R\Cent_P(T))$.

Since $R\Cent_P(T)$ is fully $\F$-normalized and $K$ is a normal $p$-subgroup of $\Aut_\F(R\Cent_P(T))$, $K\leq \Aut_P(R\Cent_P(T))$, so that $c_g^\psi c_h^{-1}$, which acts trivially on $R$ and $Y_R$, lies in $\Aut_P(R\Cent_P(T))$. Therefore, $c_g^\psi\in\Aut_P(R\Cent_P(T))$, and so $N_\psi$ contains $c_g$. Thus $\Norm_T(Q)\leq N_\psi$, so as $R\Cent_P(T)$ is fully $\F$-normalized, $\psi$ extends to $\chi:N_\psi\to P$ in $\F$, which restricts to $\theta:\Norm_T(Q)\Cent_P(T)\to P$. As $Q\theta=R$, $\Norm_T(Q)\theta\leq \Norm_T(R)$, so Proposition \ref{prop:extnstoQCP(T)}(i) implies that the image of $\theta$ is contained in $\Norm_T(R)\Cent_P(T)$, as claimed.
\end{pf}

\begin{cor}\label{cor:canextendiso} Let $\phi:Q\to R$ be an $\F$-isomorphism, where $Q$ and $R$ are $\E$-centric subgroups of $T$. Let $\psi$ be an extension of $\phi$ to $Q\Cent_P(T)$ in $\F$. Suppose that $\Aut_T(Q)^\phi=\Aut_T(R)$, so that $|\Norm_T(Q)|=|\Norm_T(R)|$ (this happens for example if $\phi$ is the restriction of an automorphism of $T$). There is an extension $\theta$ of $\psi$ in $\F$, with $\theta:\Norm_T(Q)\Cent_P(T)\to \Norm_T(R)\Cent_P(T)$ an isomorphism.
\end{cor}
\begin{pf} If $X$ is a fully $\F$-normalized subgroup conjugate to $R\Cent_P(T)$, then $X=S\Cent_P(T)$ for $S=T\cap X$, by Proposition \ref{prop:extnstoQCP(T)}(i). Let $\chi:\Norm_T(R)\Cent_P(T)\to\Norm_T(S)\Cent_P(T)$ be a morphism in $\F$, with $R\chi=S$. (Such a morphism exists since there is an isomorphism $\alpha:R\Cent_P(T)\to S\Cent_P(T)$ in $\F$ such that $N_\alpha=\Norm_P(R\Cent_P(T))$, which clearly contains $\Norm_T(R)\Cent_P(T)$; the image of $\chi$ is clear by Proposition \ref{prop:extnstoQCP(T)}(i) again.)

Consider the maps $\phi\chi:Q\to S$ and $\psi\chi:Q\Cent_P(T)\to S\Cent_P(T)$. Since $N_{\phi\chi}$ contains $\Norm_T(Q)$, so does $N_{\psi\chi}$ by Proposition \ref{prop:caseRfullnormal}, and so there is an extension $\theta_1:\Norm_T(Q)\Cent_P(T)\to \Norm_T(S)\Cent_P(T)$ of $\psi\chi$ in $\F$. We claim that $\Norm_T(Q)\theta_1=\Norm_T(R)\chi$. If this is true, then $\theta=\theta_1\chi^{-1}:\Norm_T(Q)\Cent_P(T)\to \Norm_T(R)\Cent_P(T)$ is an isomorphism (since the two groups have the same order) that extends $\psi$ (since on $Q\Cent_P(T)$ it acts as $\psi\chi\chi^{-1}=\psi$).

Write $N=\Norm_T(Q)\theta_1$ and $M=\Norm_T(R)\chi$. The subgroups $N$ and $M$ are determined by the maps $\theta_1|_Q$ and $\chi|_R$ respectively, since $N$ and $M$ are the preimages in $\Aut_T(S)$ of the subgroups $\Aut_T(Q)^{\theta_1|_Q}$ and $\Aut_T(R)^{\chi|_R}$. However, since $\Aut_T(Q)^\phi=\Aut_T(R)$, we see that 
\[ \Aut_T(Q)^{\theta_1}=\Aut_T(Q)^{\phi\chi}=\Aut_T(R)^\chi,\]
so that $N=M$.
\end{pf}

We will use this corollary to prove that $\Norm_T(Q)$ is a detecting subgroup whenever $Q$ is, at least for $p'$-automorphisms.

\begin{prop} Let $\alpha$ be a $p'$-automorphism in $\Aut_\E(T)$. If $Q$ is a detecting subgroup for $\alpha$, so is $\Norm_T(Q)$.
\end{prop}
\begin{pf} Since $\F$ is saturated, $\alpha$ extends to some automorphism $\beta\in\Aut_\F(T\Cent_P(T))$. Since $\beta$ acts as an automorphism of $\Cent_P(T)$ as well, by raising $\beta$ to a suitable power, we may assume that $\beta|_{\Cent_P(T)}$ is a $p'$-automorphism as well.

Let $Q$ be a detecting subgroup for $\alpha$, so that there exists an isomorphism $\phi:Q\Cent_P(T)\to R\Cent_P(T)$ in $\F$ such that $\phi|_Q=\alpha|_Q$ and $[\phi,\Cent_P(T)]\leq \centre R$. By Corollary \ref{cor:canextendiso}, there is a map $\theta:\Norm_T(Q)\Cent_P(T)\to \Norm_T(R)\Cent_P(T)$ in $\F$ extending $\phi$. Consider the isomorphism $\gamma:\Norm_T(Q)\Cent_P(T)\to \Norm_T(R)\Cent_P(T)$ that is the restriction of $\beta$ to $\Norm_T(Q)\Cent_P(T)$, and the automorphism $\chi=\theta^{-1}\gamma$ of $\Norm_T(R)\Cent_P(T)$.

Since $\theta$ and $\gamma$ are both extensions of the same morphism on $Q$, $\chi|_R$ is the identity. Therefore, since $\chi$ induces an automorphism of $\Norm_T(R)$, and $\Cent_T(R)\leq R$, by Lemma \ref{lem:extnpauto}, $\chi|_{\Norm_T(R)}$ is a $p$-automorphism. Since both $\theta$ and $\gamma$ map $Y_Q$ to $Y_R$, $\chi$ induces an automorphism of $Y_R$; we claim that this is a $p'$-automorphism. To see this, notice that for $x\in\Cent_P(T)$, we have $(\centre Qx)\theta=\centre Rx$ and $(\centre Qx)\gamma=\centre R(x\beta)$. Therefore $(\centre Rx)\chi=\centre R(x\beta)$, so that the order of the automorphism that $\chi$ induces on $Y_R$ divides the order of $\beta|_{\Cent_P(T)}$, a $p'$-automorphism.

Therefore there exists an integer $n$ such that $\chi^n|_{\Norm_T(R)}=\chi|_{\Norm_T(R)}$ and $\chi^n$ acts like the identity on $Y_R$. Consider $\theta\chi^n$: this map acts like $\gamma$ on $\Norm_T(Q)$ (since it acts like $\theta\chi^n=\theta\chi=\theta\theta^{-1}\gamma=\gamma$), and since $\chi$ acts like the identity on $Y_R$, $\theta\chi^n$ acts like $\theta$ on $Y_Q$; i.e., $[\theta\chi^n,\Cent_P(T)]\leq T$. Hence $\theta\chi^n$ is an extension of $\alpha|_{\Norm_T(Q)}$ to $\Norm_T(Q)\Cent_P(T)$ in $\F$ such that $[\theta\chi^n,\Cent_P(T)]\leq \centre{\Norm_T(Q)}$ (using Lemma \ref{lem:equivtorightaction}). Thus $\Norm_T(Q)$ is a detecting subgroup for $\alpha$, as required.
\end{pf}

From here it is easy to prove the theorem. Recall from the start of this section that $W$ is the set of all $\phi\in\Aut_\E(T)$ that extend to $\bar\phi\in\Aut_\F(T\Cent_P(T))$ with $[\bar\phi,\Cent_P(T)]\leq \centre T$. Since $\Aut_T(T)$ is a Sylow $p$-subgroup of $\Aut_\E(T)$ and $\Aut_T(T)\leq W$ by Proposition \ref{prop:extnstoQCP(T)}(iii), and all $p'$-elements of $\Aut_\E(T)$ lie in $W$ by the previous proposition, we see that $W=\Aut_\E(T)$, as needed.

\section{Weakly Normal Maps}
\label{sec:weaklynormalmaps}

In \cite{aschbacher2008}, Aschbacher defined invariant maps and normal maps, the former producing (some) invariant subsystems of a fusion system, and the latter producing every normal subsystem of a saturated fusion system. In this section we will define a `weakly normal' map, which will in fact be basically a normal map \cite[Remark 7.5]{aschbacher2008} with one axiom removed. In this section we will prove that every weakly normal map determines a weakly normal subsystem, and that every weakly normal subsystem gives rise to a unique weakly normal map, just as with normal subsystems and normal maps.

Because the original proof in \cite{aschbacher2008} uses the fact that the normal map produces a \emph{normal} subsystem, we need to find another proof.

\begin{prop}\label{prop:identifycentricsubgroups} Let $\F$ be a saturated fusion system on a finite $p$-group $P$, and let $\E$ be an $\F$-invariant subsystem on a strongly closed subgroup $T$ of $P$. A subgroup $R\leq T$ is $\E$-centric if and only if $R$ is $\F$-conjugate to a fully $\F$-normalized subgroup $S$ of $T$ with $\Cent_T(S)\leq S$.
\end{prop}
\begin{pf} This follows from \cite[Lemma 6.5(2)]{aschbacher2008}.
\end{pf}

Related to this is the following understanding of the relationship between $\E$-conjugacy classes and $\F$-conjugacy classes of subgroups of a strongly closed subgroup $T$.

\begin{lem}\label{lem:econjfconj} Let $\F$ be a saturated fusion system on a finite $p$-group $P$, and let $\E$ be an $\F$-invariant subsystem on a strongly closed subgroup $T$ of $P$. If $Q$ is a subgroup of $T$, then the $\F$-conjugacy class $\mc Q$ containing $Q$ is a disjoint union of $\E$-conjugacy classes $\mc Q_1,\dots,\mc Q_n$, and there are automorphisms $\phi_i\in \Aut_\F(T)\leq \Aut(\E)$ of $T$ such that $\mc Q_1\phi_i=\mc Q_i$.
\end{lem}
\begin{pf} By Theorem \ref{thm:normalifffrattini}, every $\F$-morphism $\phi$ between subgroups of $T$ may be written as $\phi=\alpha\beta$, with $\alpha\in\Aut_\F(T)$ and $\beta$ in $\E$, and $\Aut_\F(T)\leq \Aut(\E)$. This clearly implies the statement.
\end{pf}

Theorem \ref{thm:normalmapprelim} in fact tells us what we want for the invariant maps of \cite{aschbacher2008} to be weakly normal.

\begin{defn}Let $\F$ be a saturated fusion system on a finite $p$-group $P$, and let $T$ be a strongly $\F$-closed subgroup. A \emph{weakly normal map} on $T$ is a function $A(-)$ on the set of subgroups $U$ of $T$, with $A(U)\leq \Aut_\F(U)$, such that
\begin{enumerate}
\item if $\phi$ is an $\F$-isomorphism whose domain is $U$, then $A(U\phi)=A(U)^\phi$,
\item if $U$ is fully $\F$-normalized, then $\Aut_T(U)\leq A(U)$,
\item $\Aut_T(T)$ is a Sylow $p$-subgroup of $A(T)$,
\item if $U$ is fully $\F$-normalized, then every element of $A(U)$ extends to an element of $A(U\Cent_T(U))$, and
\item if $U$ is fully $\F$-normalized and $\Cent_T(U)\leq U$ then for any subgroup $U\leq V\leq \Norm_T(U)$, denoting by $A(U\leq V)$ the set of automorphisms of $V$ that restrict to automorphisms of $U$, the restriction map 
\[ A(U\leq V)\to \Norm_{A(U)}(\Aut_V(U)) \]
is surjective.
\end{enumerate}
\end{defn}

The idea is that if $\E$ is the subsystem generated by the morphisms in $A(U)$ for all $U\leq T$, then $\Aut_\E(U)=A(U)$.

We will make some short remarks about this definition now. The first two conditions are simply that $A(-)$ is an invariant map in the sense of \cite[Section 5]{aschbacher2008}. The third condition is obviously necessary for $\E=\gen{A(U):U\leq T}$ to be saturated, and the fifth condition is simply the surjectivity property for fully $\F$-normalized, $\E$-centric subgroups.

The fourth condition is there to make sure that $\E$ is actually generated by the automorphisms of centric subgroups: since every automorphism of a fully $\F$-normalized subgroup $U$ extends to an automorphism of the $\E$-centric subgroup $U\Cent_T(U)$, we really do get that $\E$ is generated by automorphisms of $\E$-centric subgroups, one of the conditions of Theorem \ref{thm:normalmapprelim}.

Like with normal maps \cite[Section 7]{aschbacher2008}, because of Condition (iv), we need only define a weakly normal map on fully $\F$-normalized subgroups $U$ with $\Cent_T(U)\leq U$, and then use (i) to extend to all $\F$-conjugates of such subgroups, then (iv) to extend to all fully $\F$-normalized subgroups, and then (i) again to extend to all subgroups of $T$.

Having discussed the reasons behind the definition of a weakly normal map, we now prove the claimed result.

\begin{thm}Let $\F$ be a saturated fusion system on a finite $p$-group $P$, and let $T$ be a strongly $\F$-closed subgroup of $P$. If $A(-)$ is a weakly normal map on $T$, then there exists a weakly normal subsystem $\E$ of $\F$ on $T$ with $\Aut_\E(U)=A(U)$ for all $U\leq T$.
\end{thm}
\begin{pf} Let $\E$ denote the subsystem of $\F$ generated by $A(-)$. Since $A(-)$ is an invariant map, $\E$ is $\F$-invariant by \cite[Lemma 5.5]{aschbacher2008}, and since $\Aut_\E(T)=A(T)$ for all invariant maps, $\Aut_T(T)$ is a Sylow $p$-subgroup of $\Aut_\E(T)$. Let $\mc H$ denote the set of $\E$-centric subgroups of $T$.

\medskip

\noindent\textbf{Step 1}: \textit{$\E$ is generated by $\E$-automorphisms of elements of $\mc H$.} Since $\E$ is generated by the subgroups $A(Q)$ for $Q\leq T$, it suffices to show that each element of these is a product of (the restriction of) elements of $A(U)$, for $U$ an $\E$-centric subgroup of $T$. If $Q$ is fully $\F$-normalized then this is true by Condition (iv) of being a weakly normal map. If $R$ is a subgroup of $T$ that is $\F$-conjugate to $Q$ via $\phi$, then $\phi=\alpha\beta$, with $\alpha\in\Aut(\E)$ and $\beta$ a morphism in $\E$ by Theorem \ref{thm:normalifffrattini}. Since every element of $A(Q)$ is the restriction of an automorphism of some member of $\mc H$, so is every element of $A(Q)^\alpha=A(Q\alpha)$. Finally, if $\beta$ is an automorphism of $Q\alpha$ (so $Q\alpha=R$) then we are done; else $Q\alpha$ and $R$ are $\E$-conjugate via $\beta$, and hence $\beta$ is a composition of the restriction of automorphisms of subgroups of $T$ of larger order than that of $Q$, which by induction are generated by automorphisms of elements of $\mc H$. Therefore $A(R)=A(Q)^{\beta^{-1}}$ is generated by $\E$-automorphisms of elements of $\mc H$, completing the proof.

\medskip

Let $\mc Q$ be an $\E$-conjugacy class contained in $\mc H$, and suppose that it consists of subgroups of $T$ of largest order subject to the elements $Q$ of $\mc Q$ not satisfying $A(Q)=\Aut_\E(Q)$. Since $\Aut_\E(T)=A(T)$, $\mc Q$ consists of proper subgroups of $T$. By Condition (iv), all subgroups in $\mc H$ of larger order than $|Q|$ lie in saturated $\E$-conjugacy classes. Applying an automorphism of $\E$ if necessary (and using Lemma \ref{lem:econjfconj}) we may assume that $\mc Q$ contains a fully $\F$-normalized subgroup $Q$.

\medskip

\noindent\textbf{Step 2}: \textit{$A(Q)=\Aut_\E(Q)$.} Consider the subsystem $\Norm_\E(Q)$ of $\E$, and notice that $\Aut_\E(Q)=\Aut_{\Norm_\E(Q)}(Q)$. Since $Q$ is fully $\F$-normalized, by Proposition \ref{prop:someofHgeneration2}(i) every element of $\Aut_\E(Q)$ can be written as the product of the restriction of morphisms in $\Norm_\E(Q)$ between subgroups of $\Norm_T(Q)$ strictly containing $Q$, and all $\Norm_\E(Q)$-conjugacy classes containing such subgroups of $\Norm_T(Q)$ are saturated. We first prove that each of these maps $\phi:R\to S$ may be written as the composition of automorphisms of overgroups of $Q$ in $\Norm_\E(Q)$, proceeding by induction on $n=|\Norm_T(Q):R|$, the case $n=1$ being clear.

Let $\psi:S\to \tilde S$ be an isomorphism in $\Norm_\E(Q)$ with $\tilde S$ a fully $\Norm_\E(Q)$-normalized subgroup of $\Norm_T(Q)$. If the result holds for $\phi\psi$ and $\psi$ then it holds for $\phi=\phi\psi(\psi^{-1})$, so we must show that any $\Norm_\E(Q)$-morphism whose image is fully $\Norm_\E(T)$-normalized may be written as a product of elements of $\Aut_{\Norm_\E(Q)}(U_i)$ for various $U_i\leq \Norm_T(Q)$; hence we assume that $S$ is fully $\Norm_\E(Q)$-normalized. Since the $\Norm_\E(Q)$-conjugacy class containing $R$ is saturated and $S$ is fully normalized, $S$ is fully automized and so there is a map $\theta:R\to S$ in $\Norm_\E(Q)$ that extends to a map $\bar\theta:\Norm_{\Norm_T(Q)}(R)\to \Norm_{\Norm_T(Q)}(S)$. By induction $\bar\theta$, and hence $\theta$, is the restriction of $\Norm_\E(Q)$-automorphisms of subgroups of $\Norm_T(Q)$ strictly containing $Q$, and so $\phi$ can be expressed in such a way if and only if $\theta\phi^{-1}$ can, which is an automorphism of an overgroup of $Q$, completing the proof of the claim.

Now let $\phi$ be an automorphism in $\Aut_\E(Q)$. If $\phi$ does not lie in $A(Q)$, then it must be the composition of restrictions of $\Norm_\E(Q)$-automorphisms $\psi_i$ of subgroups $R_i$ of $\Norm_T(Q)$ strictly containing $Q$ by the previous paragraphs; it suffices to check the case where $\phi$ is the restriction of a single automorphism of some overgroup $R$. By choice of $\mc Q$ we have that $A(R)=\Aut_\E(R)$. However, by Condition (v), and the fact that $\psi\in A(Q\leq R)$, we see that $\phi\in A(Q)$, completing the proof.

\medskip

\noindent\textbf{Step 3}: \textit{In any $\E$-conjugacy class of $\E$-centric subgroups of $T$, there is a fully $\E$-normalized subgroup with the surjectivity property.} If $Q$ is a fully $\F$-normalized and $\E$-centric subgroup of $T$, then Condition (v) of being a weakly normal map implies that $Q$ has the surjectivity property. If $\alpha\in\Aut(\E)$ then $Q\alpha$ also has the surjectivity property, and by Lemma \ref{lem:econjfconj} there is such a subgroup in every $\E$-conjugacy class of $\E$-centric subgroups.

\medskip

In particular, we now know that $\E$ is saturated by Theorem \ref{thm:normalmapprelim}.

\medskip

\noindent\textbf{Step 4}: \textit{$\Aut_\E(Q)=A(Q)$ for all $Q\leq T$.} If $Q$ is fully $\F$-normalized then this follows from Step 1 and Condition (iv) of being a weakly normal map. Since $\E$ is $\F$-invariant, if $\phi:Q\to R$ is an $\F$-isomorphism with $R$ a fully $\F$-normalized subgroup of $T$ then $\Aut_\E(Q)^\phi=\Aut_\E(R)$, and since $A(Q)^\phi=A(R)$ as $A(-)$ is an invariant map, we are done.
\end{pf}

As an example, we show that the condition that $\E$ be generated by automorphisms of centric subgroups is necessary.

\begin{example} Let $P$ be the group $V_4$, and let $\F$ be the fusion system of the alternating group $A_4$. Let $\E$ be the subsystem of $\F$, on $P$, given by all $\F$-maps between subgroups of order $2$, but not their extensions to $P$. Hence $\Aut_\E(Q)$ is trivial for all $Q\leq P$, and so Condition (ii) of Theorem \ref{thm:saturationequiv} is trivially satisfied. It is also clear that $\E$ is $\F$-invariant, and that $\Aut_P(P)$ is a Sylow $2$-subgroup of $\Aut_\E(P)$. However, $\E$ is obviously not saturated, since for example it does not satisfy the conclusion of Alperin's fusion theorem.

Hence there are $\F$-invariant subsystems $\E$ with $\Aut_T(T)$ fully $\E$-automized, and with all subgroups having the surjectivity property, that are not saturated.
\end{example}

\section{The Hypercentre}
\label{sec:hypercentral}
In this section we will study the hypercentre and central extensions of fusion systems. We begin with the hypercentral subgroup theorem, proving Theorem \ref{thm:hypercentral}.

\begin{thm}[Hypercentral subgroup theorem]\label{hypercentraltheorem} Let $\F$ be a saturated fusion system on a finite $p$-group $P$, and let $Q$ and $R$ be subgroups of $P$.
\begin{enumerate}
\item If $\F=P\Cent_\F(Q)$ and $\F=P\Cent_\F(R)$ then $\F=P\Cent_\F(QR)$.
\end{enumerate}
Let $X_\F$ denote the largest (strongly $\F$-closed) subgroup of $P$ such that $\F=P\Cent_\F(X_\F)$.
\begin{enumerate}
\item[(ii)] If $Q$ is a normal subgroup of $P$ contained in $X_\F$ then $Q$ is strongly $\F$-closed, and $X_{\F/Q}=X_\F/Q$ and $\Orth_p(\F)/Q=\Orth_p(\F/Q)$.
\item[(iii)] $X_\F=\central\infty\F$.
\end{enumerate}
\end{thm}
\begin{pf} As $Q,R\leq \Orth_p(\F)$, so therefore is $QR$, and hence $\F=\Norm_\F(QR)$. It remains to show that $\Aut_\F(QR)$ is a $p$-group, but if $\phi$ is a $p'$-automorphism in $\Aut_\F(QR)$ then the restrictions to both $Q$ and $R$ must be trivial, and hence $\phi=1$. Thus $\F=P\Cent_\F(QR)$, proving (i).

Let $Q\normal P$ and $Q\leq X_\F$, and let $\phi:A\to B$ be a morphism with $A\leq Q$. The map $\phi$ extends to an automorphism of $X_\F$, which must be a conjugation map $c_g$ for some $g\in P$, and $Q^g=Q$, so that $B\leq Q$. Thus $Q$ is strongly $\F$-closed.

Let $R/Q$ be a subgroup of $P/Q$ such that $\F/Q=\Norm_{\F/Q}(R/Q)$. As $R/Q$ is normal in $P/Q$, $R\normal P$ and so $\Norm_\F(R)$ is a saturated subsystem of $\F$. By Lemma \ref{normalizersame}, $\Norm_\F(R)/Q=\Norm_{\F/Q}(R/Q)=\F/R$, and by Proposition \ref{pcfqbijection}, $\Norm_\F(R)/Q=\F/Q$ implies $\F=\Norm_\F(R)$. Therefore $\Orth_p(\F)/Q=\Orth_p(\F/Q)$.

Next, suppose in addition that $\Aut_{\F/Q}(R/Q)$ is a $p$-group (i.e., $R/Q\leq X_{\F/Q}$). Since $\Aut_\F(Q)$ is a $p$-group, any $p'$-automorphism in $\Aut_\F(R)$ must act trivially on both $Q$ and $R/Q$, and is therefore trivial; thus $\F=P\Cent_\F(R)$, and so $X_\F/Q=X_{\F/Q}$, proving (ii).

Finally, we will show that $X_\F=1$ if and only if $\centre\F=1$. One direction is immediate, so suppose that $X_\F\neq 1$; let $Z$ be the subgroup $\centre P\cap X_\F$, and note that $Z$ is non-trivial since $X_\F$ is a non-trivial normal subgroup of $P$. If $Z$ is strongly $\F$-closed, then since $\F=P\Cent_\F(Z)$ (as $Z\leq X_\F$) and $P$ acts trivially on $Z$, we actually have that $\F=\Cent_\F(Z)$, so that $Z\leq \centre\F$. If $\phi:A\to B$ is any morphism in $\F$ with $A\leq Z$, then $\phi$ extends to a $P$-automorphism of $X_\F$, which must restrict to a (trivial) automorphism of $Z$. Hence $Z$ is strongly $\F$-closed, as claimed.

Induction and (ii) of this theorem imply that $\central\infty\F\leq X_\F$. If we can show that $X_{\F/\central\infty\F}=1$ then we are done by (ii) again. However, certainly $\centre{\F/\central\infty\F}=1$, and hence $X_{\F/\central\infty\F}=1$, as needed.
\end{pf}

We now want to isolate $\central i\F$, given $\central\infty\F$; this is very easy to do.

\begin{lem}\label{lem:findcentres} Let $\F$ be a saturated fusion system on a finite $p$-group $P$. Then
\[ \central i\F=\central\infty\F\cap\central iP.\]
\end{lem}
\begin{pf} Clearly $\central i\F\leq \central iP\cap\central\infty\F$. Let $x$ be an element in $\central\infty\F\cap\central iP$, and let $\phi:\gen x\to P$ be any morphism in $\F$. We need to show that $\phi$ acts trivially on $\gen x\central{i-1}\F/\central{i-1}\F$, or in other words that $[\phi,x]\in\central{i-1}\F$, for then $\central{i-1}\F x\in\centre{\F/\central{i-1}\F}$; $\phi$ extends to an automorphism $\psi$ of $\central\infty\F$ (as $\central\infty\F\leq\Orth_p(\F)$), and by Theorem \ref{hypercentraltheorem}, $\psi=c_g$ for some $g\in P$. Since $[\psi,x]=[g,x]\in\central{i-1}P$ (as $x\in\central iP$), by induction $[\psi,x]\in\central{i-1}\F$, proving that $x\in\central i\F$.
\end{pf}

In the case of a perfect fusion system, like with a perfect group, performing a central extension twice still results in a central extension.

\begin{prop} Let $\F$ be a saturated fusion system on a finite $p$-group $P$. If $\F$ is perfect then $\central 2\F=\centre\F$.
\end{prop}
\begin{pf}Write $Z=\centre\F$, let $x$ be an element of $\central 2\F$, and for $g\in P$ let $\lambda_x$ be a map defined by $g\lambda_x=[x,g]$. We claim that this is a group homomorphism $P\to Z$ that induces a morphism of fusion systems $\F\to\F_Z(Z)$. Since $\F$ is perfect, the homomorphism $\lambda_x$ must have trivial image; i.e., $g\lambda_x=1$ for all $g\in G$, so that $x\in \centre P$. Lemma \ref{lem:findcentres} will then imply that $x\in\centre\F$, as needed.

Firstly, let us prove that $\lambda_x$ is a group homomorphism: since $x\in\central 2\F$, in particular $x\in\central 2P$, so that $[x,g]\in \centre P$ for all $g\in P$. Hence
\[ [x,gh]=[x,h][x,g]^h=[x,g][x,h],\]
so that $\lambda_x$ is a group homomorphism.

Let $X$ denote the kernel of $\lambda_x$, and notice that $Z\leq X$. Let $Q$ and $R$ be subgroups of $P$ with $Q\leq X$, and let $\phi:Q\to R$ be an $\F$-isomorphism. Since $x\in \central 2\F$, we may extend $\phi:Q\to R$ to a morphism $\psi:Q(Z\gen x)\to R(Z\gen x)$ acting as an automorphism of $Z\gen x$ that acts trivially on both $Z$ and $Z\gen x/Z$. Write $z=x(x\psi^{-1})$, so that $(zx)\psi=x$. Since $[x,Q]=1$, we have
\[ [x,Q\psi]=[zx,Q]\psi\leq([z,Q]^x[x,Q])\psi=1.\]
Therefore $R\leq X$, so that $X$ is strongly $\F$-closed. Hence there is a surjective morphism $\Phi:\F\to\E$ of fusion systems with kernel $X$, where $\E$ is a fusion system on an abelian $p$-group. If we can show that any $p'$-automorphism $\alpha$ of $P$ satisfies $\alpha\Phi=\id$ then $\E=\F_{P/X}(P/X)$. Let $g$ be an element in $P$; we need to show that $g^{-1}g\alpha$ lies in $X$; i.e., that $[x,g^{-1}g\alpha]=1$. Since $\alpha$ is a $p'$-automorphism, and $\F=P\Cent_\F(Z\gen x)$, we must have that $\alpha$ acts as the identity on $Z\gen x$, which contains both $x$ and $[x,g]$. Hence (since $[x,g^{-1}]$ is central in $P$)
\[[x,g^{-1}g\alpha]=[x,g\alpha][x,g^{-1}]^{g\alpha}=([x,g]\alpha)[x,g^{-1}]^g=[x,g][x,g^{-1}]^g=[x,g^{-1}g]=1.\]
Therefore $\E$ is the fusion system $\F_A(A)$ for some abelian $p$-group $A$. Since $\F$ is perfect, we must have $A=1$, so that $X=P$. Therefore $x\in\centre P$, as needed.
\end{pf}

\begin{prop} Let $G$ be a finite group with a Sylow $p$-subgroup $P$. If $\Orth_{p'}(G)=1$ then $\central iG=\central i\F$, where $\F=\F_P(G)$.
\end{prop}
\begin{pf} By Glauberman's $Z^*$-theorem (for $p=2$ it is from \cite{glauberman1966a}, and for odd primes it follows from the classification of the finite simple groups), if $x$ is an element of order $p$ in $P$ such that $x^G\cap P=\{x\}$ then $x\in\centre G$. The condition that $x^G\cap P=\{x\}$ is the same as $x\in \centre\F$. Therefore, $\centre G=1$ if and only if $\centre \F=1$, and hence $\central\infty G=\central\infty\F$.

By Lemma \ref{lem:findcentres}, $\central i\F=\central iP\cap \central\infty\F$. We must similarly prove that $\central iG=\central iP\cap\central\infty G$; to see this, again we only have to prove that $\centre G=\centre P\cap\central\infty G$. Certainly $\centre G\leq \centre P\cap\central\infty G$, so let $x$ be an element of $\centre P\cap\central\infty G$. Any $p'$-element $g$ of $G$ must centralize $\central\infty G$, since it acts trivially on $\central iG/\central{i-1}G$ for all $i$. Since $P$ centralizes $x$ as well, $\Cent_G(x)$ contains all $p'$-elements of $G$ and a Sylow $p$-subgroup of $G$, so that $\Cent_G(x)=G$, as claimed.

Therefore $\central iG=\central iP\cap\central\infty G=\central iP\cap\central\infty\F=\central i\F$, as claimed.
\end{pf}

\section{Intersections and Products of Subsystems}
\label{sec:intersections}

We include an application of the theory of weakly normal maps, proving the existence of a subsystem like the intersection in certain situations.

\begin{thm}\label{thm:intersectionsubsystems} Let $\F$ be a saturated fusion system on a finite $p$-group $P$, and let $\E_1$ and $\E_2$ be weakly normal subsystems of $\F$, on the strongly $\F$-closed subgroups $T$ and $\bar T$ of $P$ respectively, with $T\leq \bar T$. There exists a weakly normal subsystem $\E$, contained in $\E_1\cap\E_2$ and denoted by $\E_1\curlywedge \E_2$, such that for any fully $\F$-normalized subgroup $Q$ of $T$ with $\Cent_T(Q)\leq Q$, we have
\[ \Aut_\E(Q)=\Aut_{\E_1}(Q)\cap\Aut_{\E_2}(Q).\]
\end{thm}
\begin{pf} Let $A_1(-)$ and $A_2(-)$ be the weakly normal maps corresponding to $\E_1$ and $\E_2$ respectively, and let $A(-)$ be the map given by $A(Q)=A_1(Q)\cap A_2(Q)$ for $Q$ a subgroup of $T$ with $\Cent_T(Q)\leq Q$, and then extended to all subgroups using (i) and (iv) of the definition of a weakly normal map, as described earlier in Section \ref{sec:weaklynormalmaps}. We will show that $A(-)$ is again a weakly normal map; clearly the subsystem $\E$ generated by $A(-)$ will satisfy the conclusions of the theorem. We prove each of the five conditions in turn.

The first and second properties are clear from the fact that each $A_i(-)$ is a weakly normal map. Clearly, since $\Aut_T(T)$ is a Sylow $p$-subgroup of $A_1(T)$ and contained in $A_2(T)$, it is a Sylow $p$-subgroup of $A(T)$, proving the third condition. The fourth condition follows by the construction of $A(-)$, and so it remains to prove the fifth condition.

Let $Q$ be a fully $\F$-normalized subgroup of $T$ with $\Cent_T(Q)\leq Q$, and let $R$ be an overgroup of $Q$ contained in $\Norm_T(Q)$. Let $\phi$ be an automorphism in $\Norm_{A(Q)}(\Aut_R(Q))$; since both $\E_i$ are weakly normal subsystems of $\F$, there are elements $\psi_i$ in $A_i(R)$ that extend $\phi$. (The existence of $\psi_1$ is clear, and there is obviously a map $\psi_2:R\to P$ extending $\phi$. To see that $\psi_2$ is an automorphism of $R$, apply the same argument used at the end of the proof of Corollary \ref{cor:canextendiso}.)

Let $S$ be a fully $\F$-normalized subgroup of $T$ that is $\F$-conjugate to $R$ via $\alpha:T\to S$. Notice that $\Aut_T(S)^\alpha$ is a Sylow $p$-subgroup of $A_1(R)$ and is contained in $A_2(R)$. Hence $A_2(R)$ is a normal subgroup of $p'$-index of $A=A_1(R)A_2(R)$, and so all $p$-subgroups of $A$ are $p$-subgroups of $A_2(R)$.

Let $K$ be the subgroup of $A$ consisting of all automorphisms of $R$ acting trivially on $Q$; $K$ is a $p$-subgroup of $A$ by Lemma \ref{lem:extnpauto}. Since all $p$-subgroups of $A$ are $p$-subgroups of $A_2(R)$, we see that $K$ is contained in $A_2(R)$. Notice that $\psi_1\psi_2^{-1}$ acts trivially on $Q$ so lies in $K$, and therefore lies in $A_2(R)$. Hence $\psi_1\in A_2(R)$, so $\psi_1\in A(R)$, as needed. Hence Condition (v) of being a weakly normal map is satisfied.
\end{pf}

We can slightly relax the conditions of Theorem \ref{thm:intersectionsubsystems} to allow the case where $\E_1$ is just a saturated subsystem, and not $\F$-invariant. In this case, $\E$ is weakly normal in $\E_1$, but not necessarily weakly normal in $\F$.

We get two corollaries of Theorem \ref{thm:intersectionsubsystems}.

\begin{cor} Let $\F$ be a saturated fusion system on a finite $p$-group $P$, and let $\E_1$ and $\E_2$ be weakly normal subsystems of $\F$, on the same strongly $\F$-closed subgroup $T$ of $P$. Then $\Orth^{p'}(\E_1)=\Orth^{p'}(\E_2)$.
\end{cor}

This shows that the set of weakly normal subsystems on a given subgroup $T$, partially ordered by inclusion, has a unique minimal element, which we will denote by $\mc R_\F(T)$. We also get another corollary, relating these minimal weakly normal subsystems.

\begin{cor} Let $\F$ be a saturated fusion system on a finite $p$-group $P$, and let $T$ and $\bar T$ be strongly $\F$-closed subgroups of $P$, with $T\leq \bar T$. If there are weakly normal subsystems of $\F$ on both $T$ and $\bar T$, then $\mc R_\F(T)\leq \mc R_\F(\bar T)$.
\end{cor}

In fact, there is also a largest weakly normal subsystem on a given strongly closed subgroup, if the set of such weakly normal subsystems is non-empty. This follows from a result of Puig.

\begin{thm}[Puig {{\cite[Proposition 6.10]{puig2006}}}]\label{thm:chooseweaklynormal} Let $\F$ be a saturated fusion system on a finite $p$-group $P$, and let $\E$ be a weakly normal subsystem of $\F$, on a strongly $\F$-closed subgroup $T$ of $P$. If $H$ is any normal subgroup of $\Aut_\F(T)$ containing $\Aut_\E(T)$, with $|H:\Aut_\E(T)|$ prime to $p$, then there exists a weakly normal subsystem $\E'$ of $\F$ on $T$, containing $\E$, and such that $\Aut_{\E'}(T)=H$. Furthermore, for any $Q\leq T$,
\[ \Hom_{\E'}(Q,T)=\Hom_\E(Q,T)\cdot H.\]
\end{thm}

In particular, if we take $\E=\mc R_\F(T)$, and $H$ to be the largest such subgroup of $\Aut_\F(T)$, then we get the largest weakly normal subsystem of $\F$ on $T$. We will denote this by $\mc R^\F(T)$. (Notice that if $\E'$ is any weakly normal subsystem of $\F$ on $T$, then $\E'$ contains $\E$ and $|\Aut_{\E'}(T):\Aut_\E(T)|$ is prime to $p$, so every weakly normal subsystem of $\F$ on $T$ can be obtained by the method in Theorem \ref{thm:chooseweaklynormal}.)

If there are (weakly) normal subsystems of $\F$ on $T$, then we say that $T$ is \emph{based} in $\F$. (In the next section we note that in \cite{aschbacher2008}, Aschbacher proves that there are strongly $\F$-closed subgroups that are not based.) In \cite{aschbacher2007un}, Aschbacher proves the following theorem.

\begin{thm}[Aschbacher] Let $\F$ be a saturated fusion system on a finite $p$-group $P$, and let $T_1$ and $T_2$ be strongly $\F$-closed subgroups. If $T_1$ and $T_2$ are based, so is $T_1\cap T_2$.
\end{thm}

An obvious question is to ask whether, in this case, $T_1T_2$ is also based. In the same article, Aschbacher proves that this is true if $T_1$ and $T_2$ commute, but the general case is still open. We will use the above theorem to develop a general theory of intersections of weakly normal subsystems.

Let $T$ be a strongly $\F$-closed subgroup of $P$ that is based in $\F$, and let $\E$ be a subsystem (not necessarily saturated) of $\F$ on some subgroup $Q$ of $P$, with $T\leq Q$. Suppose that $\E$ contains $\mc R_\F(T)$. It is easy to see via Theorem \ref{thm:chooseweaklynormal} that the set of all weakly normal subsystems of $\F$ on $T$ that are contained in $\E$ has a unique largest element, analogous to the core of a subgroup of a finite group; we will call such a subsystem the \emph{$T$-core of $\E$}. If $\E_1$ and $\E_2$ are any two subsystems of $\F$ on subgroups $Q_1$ and $Q_2$ of $P$, such that $T=Q_1\cap Q_2$ is a strongly $\F$-closed subgroup of $P$ based in $\F$, then we can construct the unique maximal weakly normal subsystem of $\F$ contained in both $\E_i$, simply the $T$-core of $\E_1\cap\E_2$, as long as $\mc R_\F(T)\leq \E_1\cap\E_2$. In particular, if both $\E_i$ are weakly normal then such a subsystem exists, and as above we will denote it by $\E_1\curlywedge\E_2$. (This notation extends the earlier one, as when both `intersections' are definable they coincide.) This yields a theory of intersections of weakly normal subsystems.

It should be noted that if the $\E_i$ are both normal in $\F$, $\E_1\curlywedge\E_2$ need not be normal in $\F$. Indeed, the intersection $\E_1\wedge\E_2$ constructed by Aschbacher in \cite{aschbacher2007un} \emph{does not} coincide with the subsystem constructed here in general. We will see an example of such a situation in the next section.
\bigskip

Having proved that if $T\leq \bar T$ are strongly $\F$-closed subgroups that are based in $\F$ then $\mc R_\F(T)\leq \mc R_\F(\bar T)$, we turn our attention to a related question, whether $\mc R^\F(T)\leq \mc R^\F(\bar T)$. For a satisfactory theory of products of weakly normal subsystems to be constructed, this must be true, since else there need not be any weakly normal subsystem on $\bar T$ containing both $\mc R^\F(T)$ and $\mc R_\F(\bar T)$, an obvious necessity to define a product of those two subsystems.

However, as the following example shows, it is not true in general that $\mc R^\F(T)\leq \mc R^\F(\bar T)$.

\begin{example} Let $Q$ be elementary abelian of order $9$, generated by $a$ and $b$. Let $\phi$ be the map fixing $b$ and sending $a$ to $ab$, and let $\psi$ be the map fixing $a$ and inverting $b$. The group $H=\gen{\psi,\phi}$ is isomorphic with $S_3$; let $G=Q\rtimes H$ and write $P$ for the (unique) Sylow $3$-subgroup of $G$. Finally, write $R=\gen b$. Notice that $R\normal G$, so that $\F_R(R)$ is a normal subsystem of $\F=\F_P(G)$. It is easy to see that the subsystem $\F_R(\gen{R,\psi})$ is isomorphic with the fusion system of $S_3$ and is equal to $\mc R^\F(R)$. However, since $\Aut_\F(Q)\cong S_3$, there can be no weakly normal subsystems of $\F$ on $Q$ except for $\F_Q(Q)$, as there are no normal $3'$-subgroups of $\Out_\F(Q)=\Aut_\F(Q)$. Thus $\mc R_\F(Q)=\mc R^\F(Q)=\F_Q(Q)$. Therefore
\[ \mc R^\F(R)\not\leq \mc R^\F(Q).\]
(There are three saturated subsystems of $\F$ that contain both $\mc R^\F(R)$ and $\mc R_\F(Q)$, namely those given by the subgroups $\gen{Q,\psi\phi^i}$ for $i=0,1,2$. However none of these is weakly normal in $\F$.)
\end{example}

This means that there can be no `reasonable' theory of products of weakly normal subsystems of fusion systems. It might still be possible to produce a theory of products of \emph{normal} subsystems of fusion systems.

\bigskip

We end with a couple of examples of situations that show that some obvious candidates for an `intersection subsystem' do not work. To see why we cannot simply take $\E_1\cap\E_2$, consider the following example.

\begin{example} Let $G$ be the group $A_4\times A_4$, and let $P$ be the Sylow $2$-subgroup of $G$, an elementary abelian group of order $16$, generated by $a$ and $b$ in the first factor, and $c$ and $d$ in the second. Let $x$ and $y$ be elements of order $3$, with both sending $a$ to $b$, and $x$ sending $c$ to $d$ and $y$ sending $d$ to $c$. Let $H_1$ be the group isomorphic with $A_4$ generated by $P$ and $x$, and $H_2$ be the group generated by $P$ and $y$.

Let $\F$ be the fusion system of $G$ on $P$, and $\E_i$ be the subsystem generated by $H_i$. Since $H_i\normal G$, the $\E_i$ are weakly normal in $\F$, but if $\E=\E_1\cap\E_2$, then $\Aut_\E(P)$ is trivial, but there is a non-trivial automorphism on $Q=\gen{a,b}$, which therefore cannot extend to $P$. Hence $\E$ is not saturated.
\end{example}

This example shows why we should define the weakly normal map in Theorem \ref{thm:intersectionsubsystems} only on the fully $\F$-normalized, $\F$-centric subgroups, and then extend it in the unique way, rather than try to define it on all subgroups to begin with, as the subgroup $Q$ above was fully $\F$-normalized, and both $\Aut_{\E_i}(Q)$ were the same, but the `correct' choice for $\Aut_\E(Q)$ was the trivial group.

If the subsystems $\E_i$ do not lie on the same subgroup then the construction in Theorem \ref{thm:intersectionsubsystems} does not work well; firstly because the $\E$-centric subgroup $R$ need not be $\E_i$-centric, but even in this case things can go wrong, even in a fusion system $\F_P(P)$.

\begin{example} Let $P=D_8\times C_2$, with the $D_8$ factor generated by an element $x$ of order $4$ and $y$ of order $2$, and the $C_2$ factor being generated by $z$. Let $Q=\gen{x,y}$, and $R=\gen{xz,y}$; then $S=Q\cap R=\gen{x^2,y}$ is a normal Klein four subgroup of $P$. If $\F=\F_P(P)$, and $\E_1$ and $\E_2$ are the subsystems $\F_Q(Q)$ and $\F_R(R)$ respectively, then both $\E_1$ and $\E_2$ are weakly normal in $\F$. We see that $\Aut_{\E_1}(S)=\Aut_{\E_2}(S)$. However, the only saturated subsystem of $\F$ on $S$ is $\E=\F_S(S)$, for which $\Aut_\E(S)$ is trivial, and so the `correct' subsystem we want inside $\E_1\cap \E_2$ is $\E$ itself. Taking the intersection of the $\Aut_{\E_i}(S)$ therefore does not yield a saturated subsystem.

(Notice that in this case, $S$ is both fully $\F$-normalized (indeed, strongly $\F$-closed), and $\E_i$-centric for $i=1,2$, so the problem does not lie in not being $\E_i$-centric.)
\end{example}

\section{Corollaries of Theorem \ref{thm:weaknormalnormal}}
\label{sec:cortothma}

Here we collect a few corollaries of Theorem \ref{thm:weaknormalnormal}. We begin with the corollary mentioned in the introduction.

\begin{cor} A saturated fusion system has no proper, non-trivial weakly normal subsystems if and only if it has no proper, non-trivial normal subsystems.
\end{cor}

Using this corollary, it is easy to show that a fusion system being quasisimple (i.e., a perfect fusion system such that $\F/\centre\F$ is simple) is not dependent on which definition of simplicity is used. In \cite{aschbacher2007un}, Aschbacher proves the existence of the analogue of the generalized Fitting subgroup for fusion systems. This is the central product of $\Orth_p(\F)$ and all subnormal quasisimple subsystems. It is easy to see that any weakly subnormal (i.e., the transitive closure of being weakly normal) quasisimple subsystem is also subnormal, and therefore we have the following corollary.

\begin{cor} The generalized Fitting subsystem of a fusion system $\F$ is the central product of $\Orth_p(\F)$ and all weakly subnormal, quasisimple, subsystems.
\end{cor}

In other words, the definition of the generalized Fitting subsystem does not depend on the definition of normality used, just as the definition of a simple fusion system.

\medskip

A \emph{minimal (weakly) normal subsystem} of a fusion system $\F$ is a (weakly) normal subsystem $\E$ of $\F$ such that if $\E'$ is a (weakly) normal subsystem of $\F$ contained in $\E$, then either $\E'=\E$ or $\E'=1$. The same methods used for minimal normal subgroups prove that minimal normal subsystems are isomorphic to direct products of isomorphic simple fusion systems.

\begin{cor} A subsystem $\E$ of a saturated fusion system $\F$ is a minimal normal subsystem if and only if it is a minimal weakly normal subsystem.
\end{cor}

If $\E$ is a normal subsystem of a saturated fusion system $\F$ and $\E$ is contained in a saturated subsystem $\F'$ of $\F$, then it need not be true that $\E$ is normal in $\F'$, even in the case where $\F'$ is itself normal. The following example proves this.

\begin{example} Let $G$ be the group $S_3\times S_3$, and let $H$ denote the subgroup of index $2$ not containing either direct factor. Let $P$ denote a Sylow $p$-subgroup of $G$, and let $K$ denote the first $S_3$ factor; write $Q=K\cap P$. Since $K\normal G$, we have that $\F_Q(K)\normal\F_P(G)$. Also, as $H$ is a normal subgroup of $G$, $\F_P(H)\normal \F_P(G)$. Clearly also, $\F_Q(K)\leq \F_P(H)$. However, it is not true that $\F_Q(K)\normal\F_P(H)$. (This example also shows that the $Q$-core of $\F_Q(K)\cap \F_P(H)$, the subsystem $\F_Q(K)\curlywedge\F_P(H)$, is not the same as $\F_Q(K)\wedge \F_P(H)$, since the former is simply $\F_Q(K)$, and the latter is $\F_Q(Q)$.)
\end{example}

It is clear that always $\E$ is \emph{weakly normal} in $\F'$, and this is enough to deduce the following corollary to Theorem \ref{thm:weaknormalnormal}.

\begin{cor} Let $\F$ be a saturated fusion system on a finite $p$-group $P$, and let $\E$ be a normal subsystem of $\F$. If $\F'$ is any saturated subsystem of $\F$ containing $\E$, and $\Orth^{p'}(\E)=\E$, then $\E\normal\F'$.
\end{cor}

In \cite[Section 7]{aschbacher2008}, Aschbacher gave examples of strongly $\F$-closed subgroups on which there are no normal subsystems of $\F$ defined. In view of Theorem \ref{thm:weaknormalnormal}, we therefore have the following corollary.

\begin{cor} If $\F$ is a saturated fusion system on a finite $p$-group $P$, and $T$ is a strongly $\F$-closed subgroup of $P$, then there need not be a weakly normal subsystem of $\F$ defined on $T$.
\end{cor}

\section{Comparing Weakly Normal and Normal Subsystems}

In this short section we will give an overview of the differences between weakly normal and normal subsystems. We begin with the following easy lemma, needed to start to understand the direct products of fusion systems. (For a definition of the direct product of two fusion systems, see \cite[Section 1]{blo2003}.)

\begin{lem} Let $\F$ be a saturated fusion system on the finite $p$-group $P$, and suppose that $P$ is the direct product of $P_1$ and $P_2$, two strongly $\F$-closed subgroups of $P$. Write $\E_i$ for the full subcategory of $\F$ on $P_i$. We have that $\E_i$ is a weakly normal subsystem of $\F$, that $\F\leq \E_1\times \E_2$, and that $\F=\E_1\times\E_2$ if and only if both $\E_1$ and $\E_2$ are normal in $\F$. This last condition holds if and only if at least one of the $\E_i$ is normal in $\F$.
\end{lem}
\begin{pf} Let $\E_i$ be defined as above; obviously it suffices to prove that $\E_1$ is weakly normal in $\F$ to prove that both $\E_i$ are, so let $A(U)=\Aut_\F(U)$ for $U\leq T=P_1$. We prove that $A(-)$ is a weakly normal map: clearly (i) and (ii) are satisfied, and (iii) is satisfied since $\Aut_T(T)=\Aut_P(T)$. If $U$ is fully $\F$-normalized and $\phi\in\Aut_\F(U)$, then $\phi$ extends to an automorphism $\bar\phi$ of $U\Cent_P(U)=U\Cent_T(U)\times P_2$ in $\F$; as $T$ is strongly $\F$-closed, $\bar\phi$ restricts to an automorphism in $\Aut_\F(U\Cent_T(U))$, proving (iv).

To prove the fifth property, let $U$ be fully $\F$-normalized and $U\leq V\leq \Norm_T(U)$, and let $\phi\in\Norm_{\Aut_\F(U)}(\Aut_V(U))$. Clearly $\phi$ extends to $\bar\phi$ in $\Aut_\F(U\times P_2)$ acting trivially on $P_2$, and since $V$ acts trivially on $P_2$, we see that $\Aut_{V\times P_2}(U\times P_2)$ is normalized by $\bar\phi$. Thus we get an extension $\psi\in\Aut_\F(V\times P_2)$ of $\bar\phi$, which restricts to an $\F$-automorphism of $V$ extending $\phi$, completing the proof of (v). Thus the $\E_i$ are weakly normal in $\F$.

If both $\E_i$ are normal in $\F$ then $\E_1\times\E_2$ is a normal subsystem of $\F$ by \cite[Theorem 3]{aschbacher2007un}, and by the definition of the $\E_i$, we see that $\F=\E_1\times \E_2$. If $\F=\E_1\times \E_2$ on the other hand, it is easy to see that every automorphism of each $P_i$ extends to an automorphism of $P$ acting trivially on $P_{3-i}$, so this equivalence is proved. Finally, suppose that $\E_1$ is normal in $\F$, and let $\phi_2\in\Aut_\F(P_2)$. Since $\F$ is saturated, it extends to some automorphism $\phi=(\phi_1,\phi_2)$ on $P_1\times P_2$, where this notation means that $\phi_i\in\Aut(P_i)$. Since $\E_1\normal\F$, the automorphism $(\phi_1,1)$ lies in $\F$, and so $(1,\phi_2)$ also lies in $\F$, proving that $\E_2\normal\F$.
\end{pf}

Hence if one wants direct products to work, one needs normal subsystems, and not merely weakly normal subsystems. Another deficiency in the definition of weakly normal subsystems is to do with $p$-power extensions of a fusion system. Let $\E$ be a saturated fusion system on a finite $p$-group $T$, and suppose that $P$ is a $p$-group containing a normal subgroup isomorphic to $T$. (We will identify $T$ with this subgroup.) Suppose that all elements of $\Aut_P(T)$ induce automorphisms of $\E$. In \cite[Section 4.2]{bcglo2007}, it was shown how to construct a saturated fusion system $\F'$ on $P$ containing $\E$, and such that $T$ is strongly $\F'$-closed and $\F'/T=\F_{P/T}(P/T)$. In other words, this is something like an extension of $\E$ by $P/T$.

In \cite[Theorem 5]{aschbacher2007un}, Aschbacher went further, and showed that if $\F$ is any saturated fusion system on $P$ containing $\E$ as a normal subsystem then $\F'\leq\F$. Since the only requirement for the construction of $\F'$ is that $\Aut_P(T)\leq\Aut(\E)$ (in other words, something satisfied by weak normality), one might try to relax this condition to $\E$ being weakly normal in $\F$.

\begin{prop} Let $\F$ be a saturated fusion system on a finite $p$-group $P$, and let $\E$ be a weakly normal subsystem of $\F$, on the subgroup $T$ of $P$. Let $\F'$ be the saturated fusion system on $P$ containing $\E$ such that $\F'/T=\F_{P/T}(P/T)$. We have that $\F'\leq \F$ if and only if $\E$ is normal in $\F$.
\end{prop}
\begin{pf} One direction of this proof is the result above, so assume that $\F'\leq \F$. Let $\phi\in\Aut_\E(T)$ be a $p'$-automorphism; since $\F'$ is saturated, this extends to $\bar\phi\in\Aut_{\F'}(T\Cent_P(T))$ which we may choose to be a $p'$-automorphism, and since $\F'/T=\F_{P/T}(P/T)$, $\bar\phi$ must act as an inner automorphism on $P/T$, which must be trivial since $\bar\phi$ is a $p'$-automorphism. Hence $[\bar\phi,\Cent_P(T)]\leq\centre T$. If $\phi$ is a $p$-automorphism of $T$, then $\phi=c_g$ for some $g\in T$, so that it obviously extends to $\bar\phi=c_g\in\Aut_\F(T\Cent_P(T))$ with $[\bar\phi,\Cent_P(T)]\in\centre T$. Hence $\E\normal\F$ as required.
\end{pf}

This gives two different situations in which the notion of a normal subsystem is better than the notion of a weakly normal subsystem. The chief advantage of weakly normal subsystems is that they are easier to work with in some situations, and one may use Theorem \ref{thm:weaknormalnormal} to pass between normal and weakly normal subsystems with relative ease.

\bigskip

\noindent\textbf{Acknowledgement}: I would like to thank Bob Oliver for suggesting an idea that allowed me to complete the proof of Theorem \ref{thm:weaknormalmap}. I would also like to thank Adam Glesser for listening to several iterations of the proof of Theorem \ref{thm:weaknormalnormal} until the correct version was found.

\bibliography{references}

\end{document}